\documentclass[10pt,reqno]{article}

\newif\iffigures\figurestrue
\figuresfalse

\usepackage[colorlinks=false,pdfborderstyle={/W 1}]{hyperref}
\newif\ifhyper\IfFileExists{hyperref.sty}{\hypertrue}{\hyperfalse}

\ifhyper\usepackage{hyperref}
\def\hitem#1#2{\item[\hypertarget{#1}{#2}]\expandafter\gdef\csname LBL#1ITM\endcsname{#2}}
\def\iref#1{\hyperlink{#1}{\csname LBL#1ITM\endcsname}}
\else
\def\hitem#1#2{\item[{#2}]\expandafter\gdef\csname LBL#1ITM\endcsname{#2}}
\def\iref#1{{\csname LBL#1ITM\endcsname}}
\fi

\usepackage{amssymb,amsmath,amsthm,amsfonts,mathrsfs,ulem, cancel}

\usepackage{graphicx}
\usepackage[utf8]{inputenc}

\def\text#1{\hbox{#1}}
\def\Z{{\mathbb Z}}

\def\N{{\mathbb N}}
\def\T{{\mathbb T}}

\def\P{\mathbf P}

\def\eps{\epsilon}

\def\1{\mathbf{1}}

\newtheorem{theorem}{Theorem}[section]

\newtheorem{lemma}[theorem]{Lemma}

\newtheorem*{note}{\underline{Note}}

\long\def\note#1/{\ifdraft{\marginpar{{$\Longleftarrow$}} \bf [#1] }\fi}

\newcounter{my}
\stepcounter{my}

\def\1{1\!\! 1}

\theoremstyle{definition}
\numberwithin{equation}{section}
\numberwithin{figure}{section}

\title{Invariant splitting of a slab into infinitely many robust clusters}

\author{Péter Mester\footnote{HUN-REN Alfréd Rényi Institute of Mathematics,  Budapest} 
}
\date{ \today}

\begin{document}

\maketitle

\begin{abstract}
We give an example of an invariant bond percolation process on the slab  $\Z^2\times \{0,1\}$ with the property that it has infinitely many clusters whose critical percolation probability is strictly less than $1$. We also show that no such process can exist in $\Z^2$.

\end{abstract}

\section{Introduction}

An invariant bond percolation process on a graph $G$ is a random subset of the edges whose distribution is invariant under all the graph automorphisms of $G$. 
Connected components formed by the process are called clusters, and a cluster is called {\it robust} if its critical percolation probability is strictly less than $1$ (i.e., if we delete edges of the cluster independently from each other with a small enough probability $\eps>0$, then there is still an infinite component). In this paper we show the following:

\begin{theorem}[]\label{t.slab}

There is an invariant bond process $\Phi_{\beta}({\tt RF})$ on the planar slab $\Z^2\times \{0,1\}$ which has infinitely many robust clusters.

\end{theorem}

\begin{theorem}[]\label{t.plane}
On $\Z^2$ there is no invariant process which has more than two robust clusters.

\end{theorem}

With extra conditions, such as positive associations, even having
more than one robust clusters is impossible on $\Z^2$ — see \cite{GKR} and \cite{Sh}. But our general negative result seems to be new.
Thus the slab $\Z^2\times \{0,1\}$ is radically different from the planar lattice $\Z^2$ regarding what kind of percolation processes it can support. We think this is interesting as there are papers, such as \cite{DST15a} (which shows the nonexistence of infinite clusters at critical Bernoulli percolation in the case of two-dimensional slabs $\Z^2\times \{0,1,\dots,k\}$) or \cite{NTW17} (which shows that the minimal spanning forest of a quasi-planar slab $\Z^2\times \{0,\dots,k\}$ is a tree) where it is shown that a particular process behaves in the same way on quasi-planar slabs as on the planar lattice  $\Z^2$.  While the aforementioned papers deal with processes which are well-defined and natural for a wide range of graphs, our process is tied to $\Z^2\times \{0,1\}$, and we think it would be interesting to see if a ``more natural process" (at least well-defined on a wide range of graphs) can also produce infinitely many robust clusters on $\Z^2\times \{0,1\}$. Regarding the corresponding question where the base graph is non-amenable, \cite{PT} shows that the Free Uniform Spanning Forest can have radically different cluster structure on different finite extensions of a tree.

In this paper all but the last section is devoted to the description of the process $\Phi_{\beta}({\tt RF})$ witnessing the truth of Theorem~\ref{t.slab}, while the last section proves Theorem~\ref{t.plane}.

In \cite{OP}, we constructed a translation invariant site process ${\tt RF}$ (for {\it Random Forks}) with the property that it has a
robust open and a robust closed cluster (we explain in Section~\ref{s.FOLDING}, how to get an invariant one from a translation invariant one). 
We will build our process out of this ${\tt RF}$ and we will make an effort to explain ${\tt RF}$ in a self-contained way:
instead of reproducing the formalism there, we will offer illustrations and point out the delicate points of the
construction. 

\begin{figure}[htbp]
\centerline{
\includegraphics[width=0.6\textwidth]{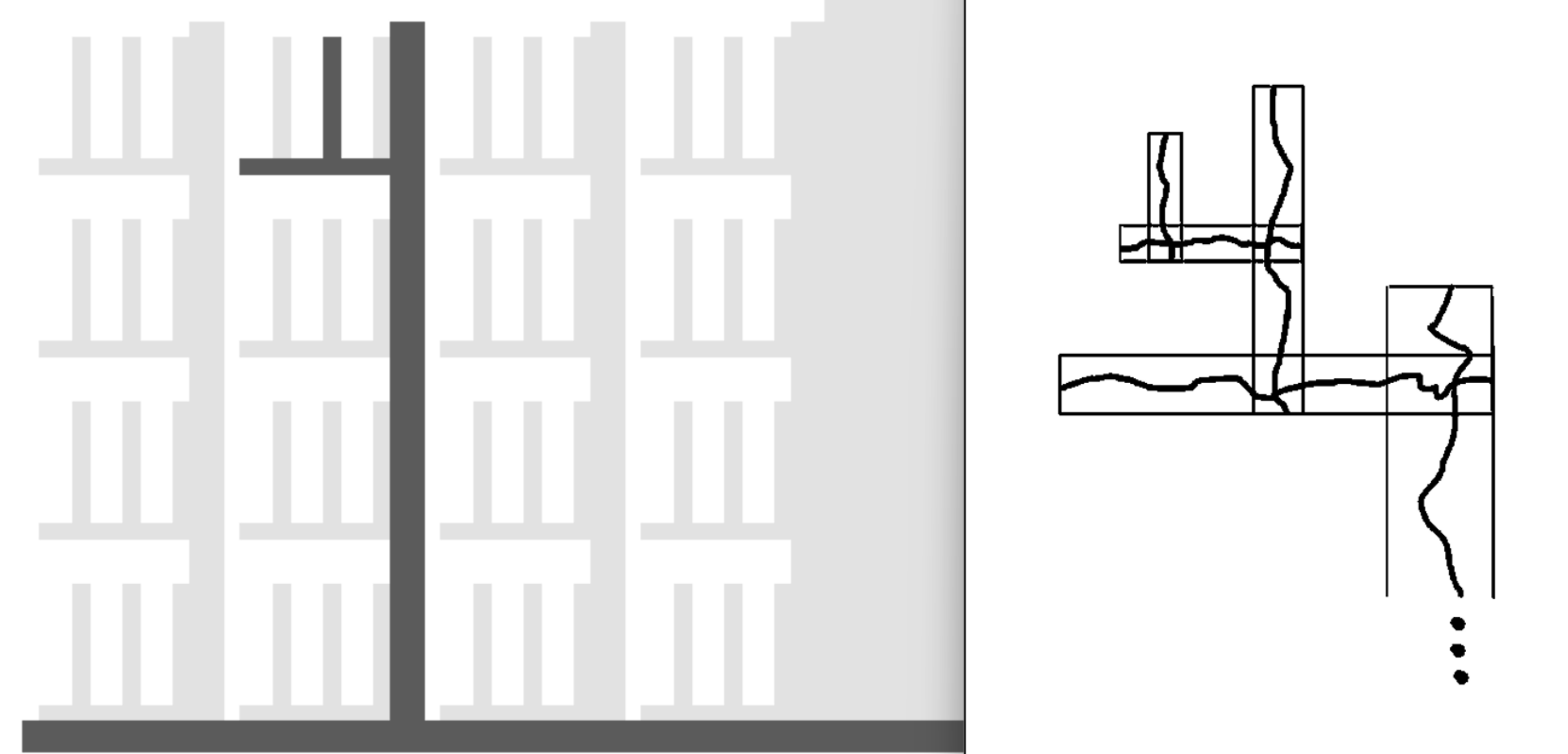}
}
\caption{The shaded part of the left picture is ${\tt RF}$, within which the dark region illustrates a ``road". The right picture illustrates that if a path exists between the opposite sides of each rectangle (in an appropriately alternating way between vertical and horizontal ones), then we can extract an infinite path from them. If the rectangles are getting ``thicker and thicker", then in a Bernoulli $p$ percolation ($p\in (1/2,1)$) this happens with positive probability due to the lower bound on the crossing probability (Lemma~\ref{l.BRcorr} below).}
\label{figChosenRoad}
\end{figure}

${\tt RF}$ is built from a collection ${\tt Rect}={\tt Rect}_{\tt V} \sqcup {\tt Rect}_{\tt H}$ (the operation is disjoint union) of vertical and horizontal rectangles which overlap in a very regular way, and from any $R\in{\tt RF}$ 
we can start a ``road" of ever widening rectangles (the intersection graph built from the rectangles of ${\tt Rect}$ is a one-ended tree, so for $R\in {\tt Rect}$ there will be a unique ${\tt next}(R)\in {\tt Rect}$ which is in the direction of infinity from $R$). In a Bernoulli $p$ percolation (for some $p>1/2$), with positive probability a so-called crossing event holds in each rectangle of a road, which implies the existence of an infinite path within the road (thus the road itself is robust). Figure~\ref{figChosenRoad} illustrates this.

\begin{figure}[htbp]
\centerline{
\includegraphics[width=0.6\textwidth]{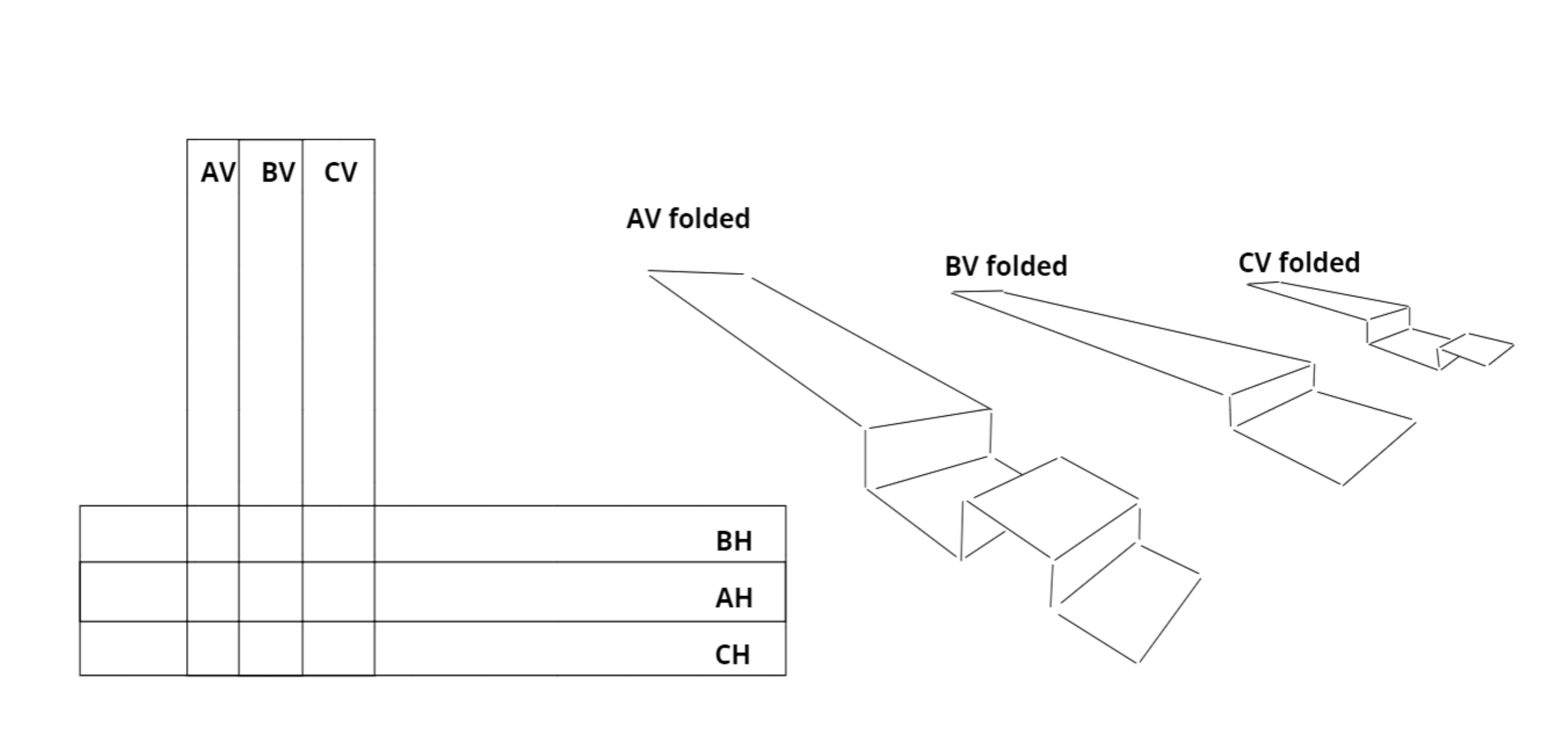}
}
\caption{Each of the three vertical slices named $AV, BV , CV$ is meant to intersect only the horizontal slice ($AH,BH, CH$) with which its first  letter agrees. In the extra space we have in the slab, we can fold each vertical slice in a way that the folded slices keep the intersections we wanted but all other intersections are avoided.}
\label{figFOLDEDslices}
\end{figure}

 In the slab $\Z^2\times \{0,1\}$ we can get infinitely many robust clusters from ${\tt RF}$ in the following way. We embed ${\tt RF}$ into the top layer of the slab, then we cut its rectangles into slices and fold those slices (using the extra space in the slab) in a way so that the folded slices only intersect specific slices of the corresponding next rectangle in the road. The result of these operation is denoted as $\Phi_{\beta}({\tt RF})$ (where $\Phi$ refers to ``folding" and $\beta$ will be a collection of maps describing the overlaps we want to maintain). By keeping a balance between not cutting the rectangles into too thin slices but still making the number of slices unbounded we can produce infinitely many robust clusters. Figure~\ref{figFOLDEDslices} illustrates schematically this folding.

\pagebreak

\section{A visual guide to ${\tt RF}$}\label{s.VISUAL}

We provide a quick visual introduction on how ${\tt RF}$ is built in \cite{OP} (most of what we cover here can be found in the section titled ``Building blocks" there).

\begin{figure}[htbp]
\centerline{
\includegraphics[width=0.3\textwidth]{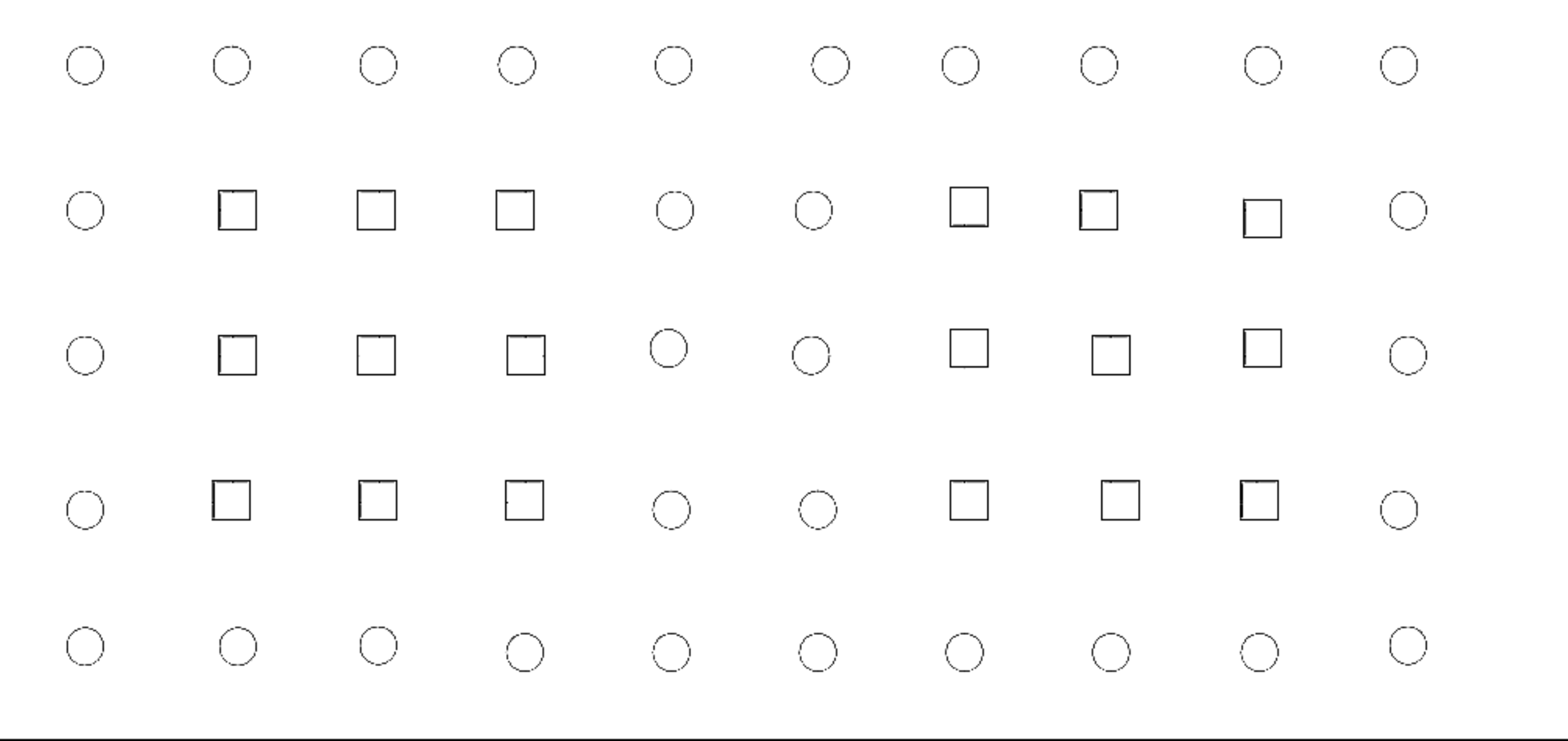}
}
\caption{A piece of a $(2,3)$-grid.}
\label{fig_2_3_grid}
\end{figure}

 A {\it grid} is a collection of congruent vertical or horizontal infinite strips of the same width which are organized periodically so that the complement of their union is a collection of congruent squares. The width $l$ of the strips and the side length $d$ of the squares determine the grid up to translation and we call such a grid an $(l,d)$-grid. There are finitely many distinct translates of an $(l,d)$-grid and selecting one of them uniformly is a translation invariant process.  Figure~\ref{fig_2_3_grid} shows part of a $(2,3)$-grid. In this picture we represented the grid using little circles for those vertices of $\Z^2$ which are in the grid and little squares for those vertices which are not. In the rest of the pictures we do not display individual vertices of $\Z^2$, instead we use solid rectangles with different shades for clearer visibility.

\begin{figure}[htbp]
\centerline{
\includegraphics[width=0.4\textwidth]{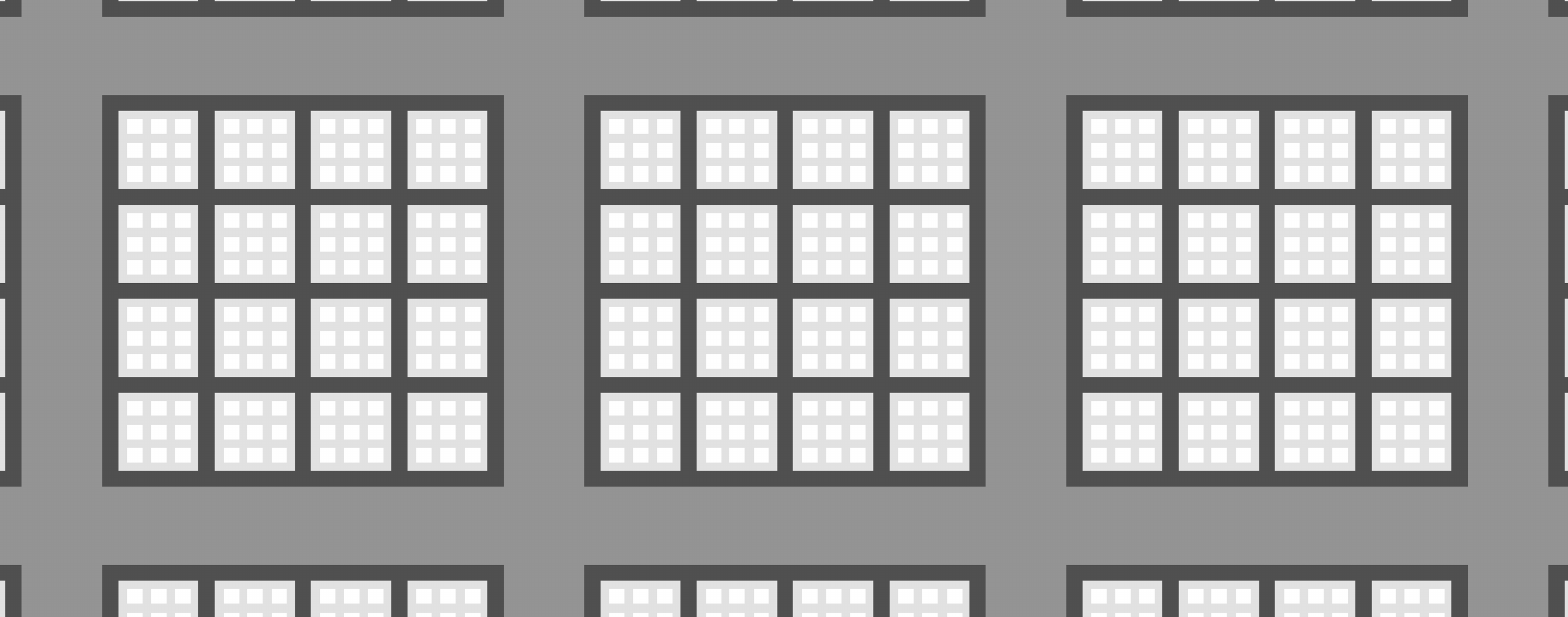}
 \hskip 1 cm
\includegraphics[width=0.3\textwidth]{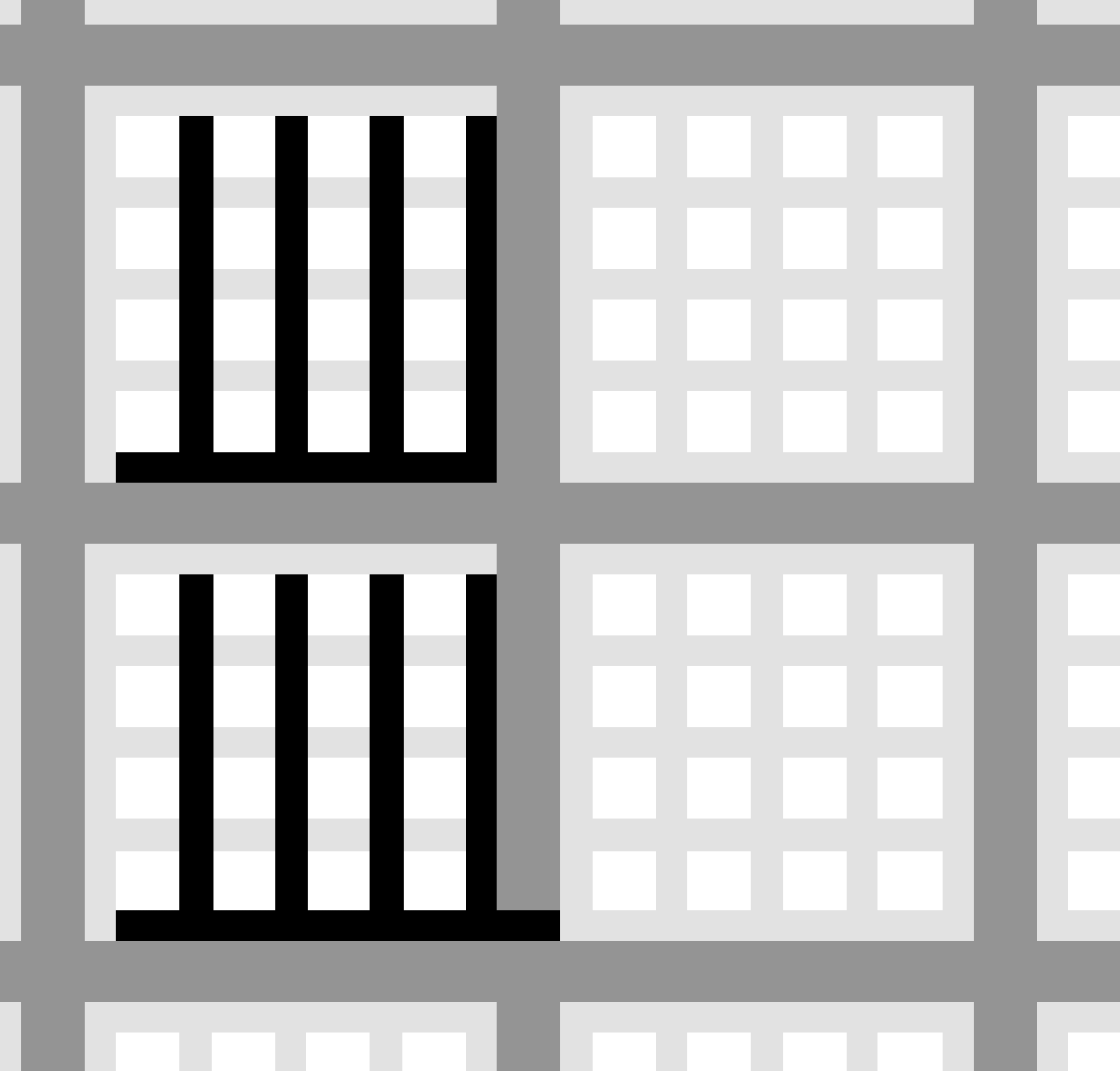}
}
\caption{The left picture shows a nested sequence of windows. The right one illustrates the fork (top left) and extended fork (bottom-left) associated to a given window in a given sequence of nested windows. The vertical and horizontal frames included in the extended forks of all the windows will constitute the sets ${\tt Rect}_{\tt V}$ and ${\tt Rect}_{\tt H}$. }
\label{figWIN}
\end{figure}

We can sequentially select grids with some fine-tuned parameters (we will see more details about those parameters in Section \ref{s.CROSSING}) so that consecutive grids fit into the previously selected collection in a way that determines a nested sequence of {\it windows}. The left part of Figure~\ref{figWIN} offers a pictorial representation of what we mean by that.

The windows are built from vertical and horizontal frames. The {\it fork} of a window is built using the all but the leftmost vertical frames and the horizontal frame at the very bottom of the window. In the top-left part of the right part of Figure~\ref{figWIN}, we highlighted the fork of a window. ${\tt RF}$ is defined by declaring all vertices/sites which belong to the fork of some window to be open. This is done for each window. Instead of selecting the fork, we could have selected the {\it extended fork} displayed at the bottom left, where the difference is that the bottom horizontal frame is extended so that it overlaps with a vertical frame from the next level window (so the extended fork is not determined by the individual window alone, but we need to know the next window as well). Even though for an individual window the fork and the extended fork are different, if we take the increasing union of them for all windows, we got the same ${\tt RF}$ process. Our reason for nevertheless introducing the extended fork (not present in
\cite{OP}) will become clear in the next section.

We can see a snapshot of ${\tt RF}$ (the shaded region of the left picture including the black rectangles) in Figure~\ref{figChosenRoad} where the windows are no longer shown. It may already suggest that the rectangles are arranged in a way which resembles a one-ended tree.

If $G$ is a graph and $X$ is a subset of its vertices, then $G[X]$ denotes the subgraph induced by $X$.
The fact that ${\tt RF}$ is a site process means that through the above process we select a subset $\bigcup_{R\in {\tt Rect}}R$ of vertices, and ${\tt RF}$ is  the subgraph induced by it; i.e., ${\tt RF}:=\Z^2[\bigcup_{R\in {\tt Rect}}R]$. In order to build the construction in this paper, we need to change this in a way we describe next.

\section{A modification of ${\tt RF}$ to define it through a graph union}\label{s.Conventions}

 Assume $G$ and $H$ are simple graphs. Then let $G\cup H$ be that simple graph whose vertex set is $V(G\cup H):=V(G)\cup V(H)$ and whose edge set is $E(G\cup H):=E(G)\cup E(H)$ (without multiplicities in the union, of course).
Sometimes we will refer to this as the ``graph union" of graphs.

By a {\it block} we mean a finite consecutive subset of integers. 
If $Q=A \times B$ where $A$ and $B$ are blocks, then we say that $Q$ is a {\it rectangle}, and we will use the notation ${\tt h}(Q)=A, {\tt v}(Q)=B$ to express the blocks in terms of the rectangle, where ${\tt h}$ and ${\tt v}$ refer to horizontal and vertical. 

If, for rectangles $Q$ and $R$, it holds that ${\tt h}(Q)\subset {\tt h}(R)$ and ${\tt v}(Q)\supset {\tt v}(R)$, and both containments are proper, then we say that the pair $(Q,R)$ is well-joined in a {\it vertical to horizontal} way. If, on the other hand, we have  ${\tt v}(Q)\subset {\tt v}(R)$ and ${\tt h}(Q)\supset {\tt h}(R)$, and both containments are proper, then we say that the pair $(Q,R)$ is well-joined in a {\it horizontal to vertical} way.

The following is a quick illustration of why we want to avoid using a site percolation.
Consider the rectangles $V_1=\{(0,0),(0,1)\},V_2=\{(1,0),(1,1)\},H_1=\{(0,0),(1,0)\}$.
 
If $G$ is the graph union of their induced subgraphs, that is $G=\Z^2(V_1)\cup \Z^2(V_2)\cup \Z^2( H_1)$, then any path in $G$ which connects pairs of vertices from $V_1$ and $V_2$ must go through  $H_1$ as well. However, if we take those rectangles just as vertex sets and declare their vertices to be open (as in a site percolation), then we got $\Z^2[V_1\cup V_2 \cup H_1]$, in which there is an edge between $(0,1)\in V_1$ and $(1,1)\in V_2$, so a path exists between them which does not go through $H_1$. In our construction later, this kind of situation would be undesirable, as we want infinitely many distinct clusters, and for that purpose it is not good if vertex disjoint rectangles are automatically connected if they are next to each other.

For this reason, the random fork subgraph in the present paper will be defined as ${\tt RF}:=\bigcup_{R \in {\tt Rect}}\Z^2[R]$ (instead of the original $\Z^2[\bigcup_{R\in {\tt Rect}}R]$).

\section{One-ended trees}\label{s.Tree}

Now we will give an abstract description of how we associate a one-ended tree to ${\tt RF}$. A tree ${\cal T}$ is one-ended if for any vertex $v$ of ${\cal T}$ there is a unique infinite simple path starting from $v$. 
In general, if we have a forest ${\cal F}$ with the property that every component of it is a one-ended tree, and $v\in V({\cal F})$ is a vertex of it, then there is a unique one--way infinite path in ${\cal F}$ which starts in $v$; we will denote by ${\tt next}_{\cal F}(v)$ the vertex right after $v$ in this path. Unless there is some ambiguity which forest we are talking about, we will omit the subscript, so we will just use the notation ${\tt next}(v)$ (note that this definition allows ${\cal F}$ to be a one-ended tree itself). 
 The rectangles in ${\tt Rect}$ will be the ones included in the extended forks of the nested window sequence which we have not yet described formally but illustrated them in the right part of Figure~\ref{figWIN}.

${\tt RF}$ can be described as a graph union of the rectangles contained in a certain collection ${\tt Rect}={\tt Rect}_{\tt V} \sqcup {\tt Rect}_{\tt H}$ (the operation is disjoint union) of rectangles, where the rectangles in ${\tt Rect}_{\tt V}$ (respectively in ${\tt Rect}_{\tt H}$) will be called {\it vertical} (respectively {\it horizontal}). 
The way ${\tt Rect}={\tt Rect}_{\tt V} \sqcup {\tt Rect}_{\tt H}$ is built makes sure that:
 Any two distinct rectangles $S_1,S_2\in {\tt Rect}_{\tt V}$ (respectively $\in {\tt Rect}_{\tt H}$) are vertex disjoint. If $S_V\in {\tt Rect}_{\tt V}$ and $S_H\in {\tt Rect}_{\tt H}$, they are either vertex disjoint or the pair $(S_V,S_H)$ is well-joined in a vertical to horizontal way or the pair $(S_H,S_V)$ is well-joined in a horizontal to vertical way. 

Let $\T_{\tt RF}$ be the graph whose vertex set is $V(\T_{\tt RF})={\tt Rect}$, and there is an edge between $S_1,S_2\in {\tt Rect}$ exactly if $S_1,S_2$ are not vertex-disjoint.

Due to the structure of the nested windows whose frames are used in the definition of ${\tt Rect}={\tt Rect}_{\tt V} \sqcup {\tt Rect}_{\tt H}$, the rectangles in ${\tt Rect}$ are indexed by the natural numbers $0,1,\dots,n,\dots$; that is, there is an indexing map ${\tt i}:{\tt Rect}\rightarrow \N$. Moreover, if $V\in{\tt Rect}_{\tt V}$, then ${\tt i}({\tt next}(V))={\tt i}(V)$, while if $H\in {\tt Rect}_{\tt H}$, then ${\tt i}({\tt next}(H))={\tt i}(H)+1$. This index ${\tt i}$ naturally belongs to a full window, and we will derive from it another index ${\tt j}:{\tt Rect}\rightarrow \N$ by defining ${\tt j}(R):=2{\tt i}(R)$, if $R\in {\tt Rect}_{\tt V}$, and ${\tt j}(R):=2{\tt i}(R)+1$, if $R\in {\tt Rect}_{\tt V}$.

The following is also true (but would not be if we used site percolation): if $v_1,v_2$ are vertices of ${\tt RF}$ and $v_1$ is a vertex of $S_1\in {\tt Rect}$ and $v_2$ is a vertex of $S_2 \in {\tt Rect}$, then any path between $v_1$ and $v_2$ in ${\tt RF}$ must go through in all the rectangles which are contained in the unique shortest path between $S_1$ and $S_2$ in ${\T}_{\tt RF}$.

Using the tree structure, we explain at an abstract level how we make infinitely many clusters from ${\tt RF}$.

 Let $\T$ be any one-ended tree and let a function $m: V(\T)\rightarrow \N^+$ be given with the property that $m(v)\leq m({\tt next}(v))$. Associate to every $v\in V(\T)$ a set $B_v$ (a ``bag" at $v$) with $|B_v|=m(v)$, and for $v\not =u$ let $B_v\cap B_u =\emptyset$. Fix for each $v$ an injection $\beta_v:B_v\rightarrow B_{{\tt next}(v)}$. Since the bags are disjoint, we can write $\beta(s)$ instead of $\beta_v(s)$ without ambiguity. Let the forest ${\cal F}_{\beta}(\T)$ be defined as follows: let the vertex set of ${\cal F}_{\beta}(\T)$ be the disjoint union $\bigcup_{v\in V(\T)}B_v$, and let us connect $s\in B_v$ with $\beta(s)\in B_{{\tt next}(v)}$ with and edge (no other edges are added). It is clear that ${\cal F}_{\beta}(\T)$ is a forest and all of its components are one-ended trees (and ${\tt next}_{{\cal F}_{\beta}(\T)}(u)=\beta(u)$); moreover, if $s_1,s_2\in B_v, s_1\not =s_2$, then by the injectivity of the $\beta$ maps, $s_1$ and $s_2$ are in different components of ${\cal F}_{\beta}(\T)$. In particular, $m(v)$ for any $v\in V(\T)$ is a lower bound on the number of clusters in ${\cal F}_{\beta}(\T)$, and if $m(v)$ is unbounded, then ${\cal F}_{\beta}(\T)$ has infinitely many clusters. 

On the other hand, if $s_1\in B_{v_1},s_2\in B_{v_2},v_1\not=v_2, {\tt next}(v_1)={\tt next}(v_2)$, then it may happen that $\beta(s_1)=\beta(s_2)$, since in this case the two $\beta$'s are actually different, the first is $\beta_{v_1}$, the second is $\beta_{v_2}$ so their individual injectivity does not prevent sending different inputs to the same output. 

Note also that the collection of injections $\beta_v :B_v\rightarrow B_{{\tt next}(v)}$ for $v\in \T$ also determines the collection of bags $B_v$ and their sizes $|B_v|=m(v)$, so the notation ${\cal F}_{\beta}(\T)$ does refer to all ingredients in the above construction.

\section{Cutting and folding}\label{s.FOLDING}

In \cite{OP} we constructed ${\tt RF}$ as a translation invariant process (which only requires invariance with respect to translations). We can transform it into a invariant one on $\Z^2$ by randomly selecting which are the vertical and horizontal directions, and also into which direction the coordinates are increasing, and when we move to the slab we can also randomly select which layer is the top and which is the bottom (in the later discussion we will assume that these choices are made). 

To realize the abstract idea of turning $\T$ into ${\cal F}_{\beta}(\T)$ in the case when $\T=\T_{\tt RF}$ we associate a bag $B_R$ to a rectangle $R\in V(\T_{\tt RF})$ in the following way. We pick a sequence $m_0\leq m_1 \leq \dots$ and we cut $R\in {\tt Rect}$ into $m_{{\tt j}(R)}$ vertical (if $R\in {\tt Rect}_{\tt V}$) or horizontal (if $R\in {\tt Rect}_{\tt H}$) slices and the union of those slices will be the bag $B_R$. Then we define a collection of injections $\beta_R:B_R\rightarrow B_{{\tt next}(R)}$ by uniformly selecting one from the finitely many possibilities, giving ${\cal F}_{\beta}(\T_{\tt RF})$. In this way $B_R$ and $B_Q$ will be disjoint for $Q\not = R$; but this only means that if $S\in B_R$, then $S\not \in B_Q$, but if $Q={\tt next}(R)$, then $S$ intersects all slices from $B_Q$. So while the associated abstract forest ${\cal F}_{\beta}(\T_{\tt RF})$ indeed consists of infinitely many different clusters if $m_i\rightarrow \infty$, if we form the graph union ${\cal G}(C):=\bigcup_{S \in V(C)} \Z^2[S]$
for a cluster $C$ of ${\cal F}_{\beta}(\T_{\tt RF})$, then for any two clusters $C_1,C_2$ of ${\cal F}_{\beta}(\T_{\tt RF})$ the vertex and edge set of the graphs ${\cal G}(C_1)$ and ${\cal G}(C_2)$ will intersect.

\begin{figure}[htbp]
\centerline{
\includegraphics[width=0.5\textwidth]{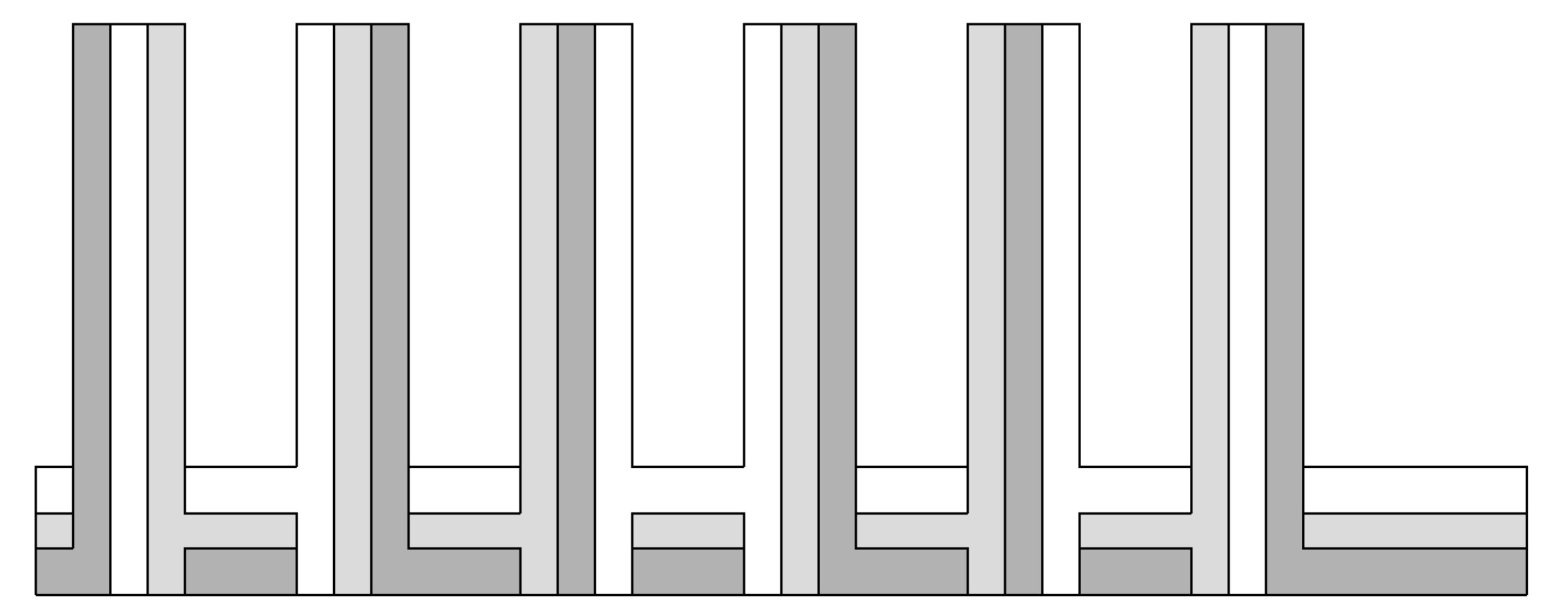}
}
\caption{Within the slab we can fold the slices (as if they were ribbons) so that only the ones of the same shade intersect.}
\label{figHarmony}
\end{figure}

Here comes the moment when embedding everything into $\Z^2 \times \{0,1\}$ helps, since we can ``fold" any slice $S\in B_R$ into a ``folded slice" $\phi(S)$ in such a way that only those folded slices have nontrivial intersection which are meant to by the $\beta$ injections. 
Figure~\ref{figHarmony} is meant to capture this idea of using some extra space in order to avoid unwanted intersections. Note that this picture is only conveying the rough idea and is not meant to be a representation of the specific folding we will use.

If $B$ is a block, $m\geq |B|$ is a positive integer, then a balanced $m$-cut of $B$ is a (possibly random) collection of $m$ blocks ${\tt cut}(m, B):=\{B_1,\dots,B_m\}$ which are pairwise disjoint whose union is $B$, and for any two $B_i,B_j\in {\tt cut}(m,B)$ we have $|B_i|-1\leq |B_j| \leq |B_i|+1$. 

This is not uniquely defined, but for our purpose any such partition will work; for example we can take the one which uniformly selects one from the finitely many possibilities.

By the horizontal $m$-cut of a rectangle $R=B_1\times B_2$ we mean the collection of rectangles (we will call them the {\it slices} of $R$):
$${\tt cut}_H(m,R):=\{{\tt h}(R) \times C: C\in {\tt cut}(m, {\tt v}(R))\}.$$
Similarly, by the vertical $m$-cut we mean the collection of rectangles (slices of $R$):
$${\tt cut}_V(m,R):=\{C\times {\tt v}(R): C\in {\tt cut}(m, {\tt h}(R))\}.$$

We will only cut a vertical rectangle into vertical slices and a horizontal rectangle into horizontal slices. 
If it is clear from the context that $R$ is a vertical rectangle, then we may just write ${\tt cut}(m,R)$ for ${\tt cut}_{\tt V}(m,R)$ (similar remarks apply for horizontal rectangles); we may also omit the subscript if we want to talk about a generic rectangle which can either be vertical or horizontal (but the cut should always be the same type as the rectangle).

Let $m_0\leq m_1\leq \dots$ be a sequence of positive integers which has the property that $m_i$ is less than the shorter side of the vertical and horizontal rectangles for the rectangles of index $i$. 

For $R\in {\tt Rect}$ let the corresponding bag be $B_R:={\tt cut}(m_{{\tt j}(R)}, R)$, then by fixing the appropriate $\beta$ injections from $B_R$ to $B_{{\tt next}(R)}$ we have the forest ${\cal F}_{\beta}(\T_{\tt RF})$.

Let us embed ${\tt RF}$ into $\Z^2\times \{0,1\}$ by placing it into the top layer $\Z^2\times \{1\}$. So if $X\subset \Z^2$, then if we use $X$ in the description of the folding operation, then we really mean $\{(x,1)\in \Z^2\times \{0,1\}:x\in X\}$.
 The folding operation will associate to any $S\in B_R$ a new connected subgraph $\phi(S)$ of $\Z^2\times \{0,1\}$ in such a way that for any two slices $S_1\in B_R,S_2\in B_Q$, for $R,Q\in {\tt Rect}$, the only way $V(\phi(S_1))\cap V(\phi(S_2)) \not = \emptyset$ if $Q={\tt next}(R)$  and $S_2=\beta(S_1)$, or the other way around ($QR={\tt next}(Q)$  and $S_1=\beta(S_2)$). Moreover, if this holds, then $V(\phi(S_1))\cap V(\phi(S_2)) = S_1\cap S_2$. This implies that if $C_1,C_2$ are two different clusters of the abstract forest ${\cal F}_{\beta}(\T_{\tt RF})$, then:

i, the graph union $\bigcup_{S\in V(C_1)} \phi(S)$ will be a connected graph (since the nonempty intersection of $S$ and $\beta(S)$ (which are neighbors in $C_1$ if $S\in V(C_1)$) is preserved by $\phi(S)\cap \phi(\beta(S))= S\cap \beta(S)$.

ii, the graph unions $\bigcup_{S\in V(C_1)} \phi(S)$ and  $\bigcup_{S\in V(C_2)} \phi(S)$ are disjoint  (and by the previous item they are connected subgraphs of $\Z^2\times \{0,1\}$).

Note that this already implies that we have infinitely many clusters if $m_i$ goes to infinity.

The formalization of the folding operation may hide its extreme simplicity so we start with an intuitive description. Given $S\in B_R$ let us place the rectangles $S$ and ${\tt next}(R)$ on the top layer of $\Z^2\times \{0,1\}$ and let us also mark the horizontal slice $\beta(S)$ as distinguished from the rest of ${\tt next}(R)$. We want to keep as much of these at the top level as we can, but we want to avoid the overlap of $S$ and that part of ${\tt next}(R)$ which does not belong to $\beta(S)$. Moreover, the slices in ${\tt next}(R)$ have a ``higher rank" than the slices in $R$ in the sense that they will stay in the top layer and the slices in $R$ will be folded in order to avoid the overlap (the slices of ${\tt next}(R)$ will be folded also in order to avoid the unwanted overlaps with the slices of ${\tt next}^2(R)$ and in this way the folding operations executed at different slices will not interfere with each other).
Imagine that ${\tt next}(R)$ is a lake but $\beta(S)$ is an island in it, and we want to build a road where $S$ lays (and in particular we want it to go through the island), but we must avoid the water and we cannot build bridges; however, we can build a subway under the lake. This subway would correspond to using the extra space in $\Z^2\times \{0,1\}$, so that we avoid the unwanted intersection by using the bottom layer.

\begin{figure}[htbp]
\centerline{
\includegraphics[width=0.4\textwidth]{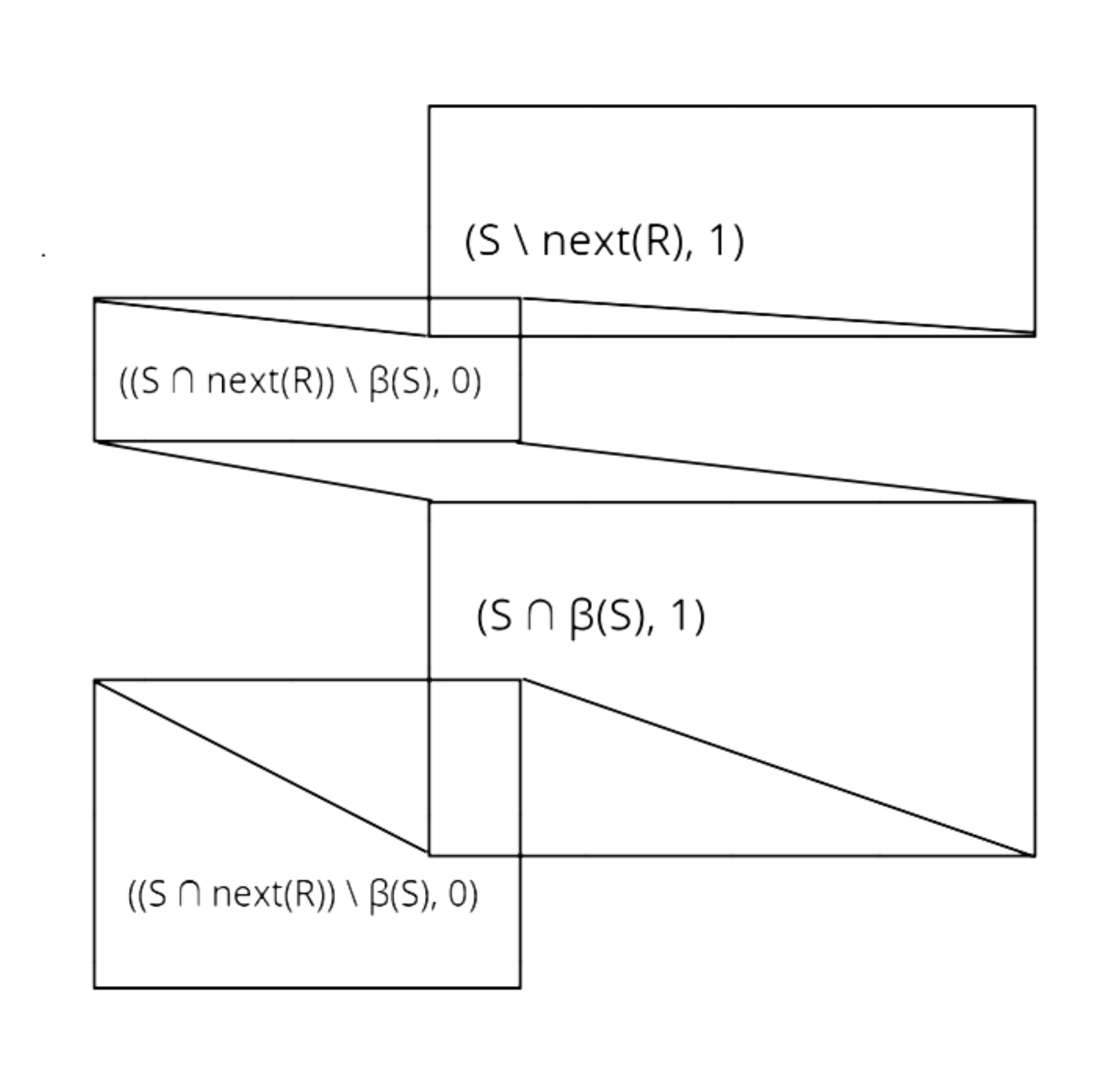}
}
\caption{A visual representation of the folding operation applied to a slice $S\in B_R$. The proportion of the pieces are highly distorted: for example, $(S\setminus {\tt next}(R),1)$ is much longer in reality than any other piece involved.}
\label{figSCHEME}
\end{figure}

Figure~\ref{figSCHEME} is meant to depict what is happening.

Before moving on, we want to point out some pathological situations which we want to avoid. If $Q\in {\tt Rect}$ and the bag $B_{{\tt next}(Q)}$ of the next rectangle ${\tt next}(Q)$ consisted of slices whose shorter side is $1$ or $2$, then in the first case the whole concept of ``folding" seems to break down while in the second case the relatively simple definition we give below would break down (the below defined sets $\partial^t$ and $\partial_b$ would force to put too many edges into the induced subgraph). So we assume that the shortest side of any slice is at least $3$, this condition can be safely assumed (the next section provide explicit details about how the parameters of the process govern the size of the rectangles in ${\tt Rect}$).

Let $S\in B_Q$, let $R=S\cap {\tt next}(Q)$. Let $I_{\phi}$ be the following indicator function $I_{\phi}:S\rightarrow \{0,1\}$: for $v\in S$, $I_{\phi}(v)= 1$ if $v\in (S\setminus R) \cup \beta(S)$ and $I_{\phi}(v)=0$ if $v\in R\setminus \beta(S)$. Let $\partial^t\subset S$ be defined by: for $v\in S$ we have $v\in \partial^t$ iff $I_{\phi}(v)=1$ but $v$ has a neighbor $u\in S$, for which $I_{\phi}(v)=0$. $\partial^t$ is meant to stand for ``top boundary" (since that part of $S$ is meant to be kept at the top layer). Let $\partial_b\subset S$ be defined by: for $v\in S$ we have $v\in \partial_b$ iff $I_{\phi}(v)=0$ but $v$ has a neighbor $u\in S$, for which $I_{\phi}(v)=1$. $\partial_b$ is meant to stand for ``bottom boundary" (since that part of $S$ is meant to be kept at the bottom layer).

Recall that if $G=(V,E)$ is a simple graph and $X\subset V(G)$, then by $G[X]$ we mean the subgraph induced by $X$. Since using this notation for $G=\Z^2\times \{0,1\}$ is quite cumbersome, we will temporarily denote $\Z^2\times \{0,1\}$ by ${\cal S}_2:=\Z^2\times \{0,1\}$ here. Let us emphasize that ``union" below means graph union.
Let $$\phi(S):={\cal S}_2[\{(v,I_{\phi}(v)):v\in S\}]\cup {\cal S}_2[\{(v,i):v\in \partial^t,i \in \{0,1\}\}]\cup {\cal S}_2[\{ (v,0):v\in \partial^t\cup\partial_b\}].$$

Now we are ready to define our construction. In order to see the full process, we go through the steps we did so far as well. 

We start with ${\tt RF}$ which depends on some sequence of parameters $(l_0,d_0),L_1,L_2\dots$ (they are positive integers and we say more about them in the next section). We add a sequence $m_0\leq m_1\leq \dots$, with $m_i\rightarrow \infty$ of positive integers as an extra parameter (which will need to grow slower than the shorter side of the rectangles in ${\tt RF}$). This defines splitting the rectangles of ${\tt RF}$ into slices, and the slices of a given rectangle can be used as its bag. 
Choosing the collection of injections to which we refer by $\beta$, the ``abstract" forest ${\cal F}_{\beta}(\T_{\tt RF})$ becomes well defined. 
Finally, let 
$$\Phi_{\beta}({\tt RF}):=\bigcup_{S\in V({\cal F}_{\beta}(\T_{\tt RF}))}\phi(S),$$ 
where the union is again meant in the graph union sense. That is, we simply replace every slice $S$ with its folded version $\phi(S)$, and then the graph whose vertices are the folded slices and where an edge represents non-empty intersection will be isomorphic to ${\cal F}_{\beta}(\T_{\tt RF})$. What we claim is that, if the parameters of ${\tt RF}$ are chosen in such a way that it is robust and its robustness can be proven by a crossing argument described in the next section, then the sequence  $m_0\leq m_1\leq \dots$ can be chosen to grow slowly enough, so that 
$\Phi_{\beta}({\tt RF})$ has infinitely many components and they are all robust.

\section{Crossing probabilities}\label{s.CROSSING}

First we say more about the parameters in ${\tt RF}$ (in \cite{OP} this is mainly covered in the Section titled ``The actual construction").

\begin{figure}[htbp]
\centerline{
\includegraphics[width=0.7\textwidth]{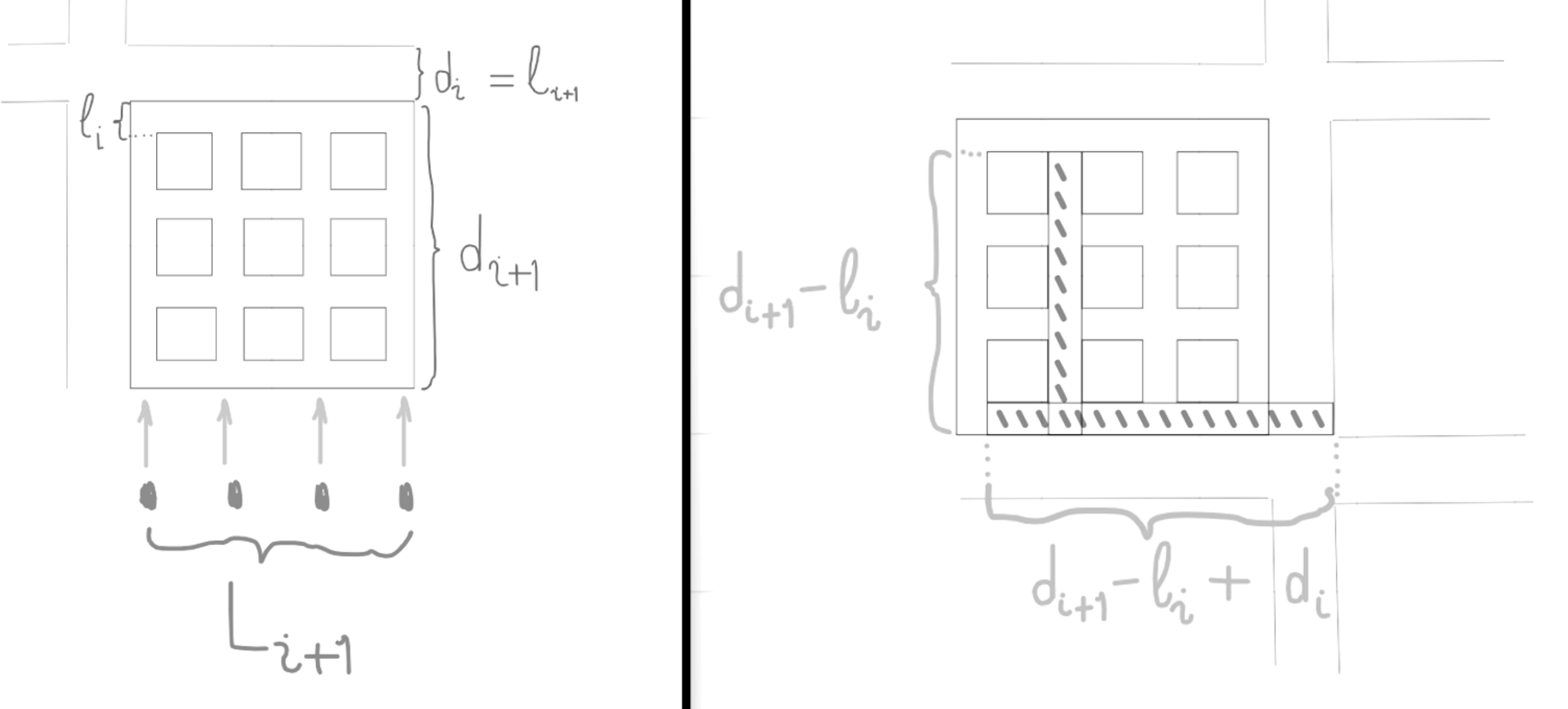}
}
\caption{From the left figure we can read off that $l_{i+1}=d_i$ and 
$d_{i+1}=L_{i+1} l_i + (L_{i+1}-1)d_i$. For the sum $l_i+d_i$ this gives the simple recursion $l_{i+1}+d_{i+1}=L_{i+1}(l_i+d_i)$. In the right picture we can see that the vertical rectangle in an extended fork is an $l_i \times (d_{i+1}-l_i)$ one while the horizontal rectangle is a $(d_{i+1}+d_i-l_i)  \times l_i$ one.}
\label{gridRECcool}
\end{figure}

To define ${\tt RF}$ we need an initial pair $(l_0,d_0)$ of positive integers and a sequence of positive integers $L_1,L_2,\dots $. The first $(l_0,d_0)$ parameter tells us (up to a random translation) what is the initial grid we start with is, then the $L_i$ sequence tells us the number of vertical (or, what is the same, horizontal) frames in the windows which will be carved out from the previous grid by the new one. This determines the parameter $(l_i,d_i)$ of the $i$th grid, as illustrated in the left part of Figure~\ref{gridRECcool}. We can also express the sizes of the vertical and horizontal rectangles of an extended fork with these parameters (see the right part of Figure~\ref{gridRECcool}). Those sizes need to be known if we want to estimate crossing probabilities.

With the right choice of the above parameters, the ${\tt RF}$ process will be well-defined and robust. To be well-defined it is important that any vertex is covered by finitely many grids only. By a Borel-Cantelli argument, the condition $\sum_{i=0}^{\infty}{1\over {L_i}}<\infty$ will be sufficient for this.  

The robustness is taken care of by a crossing argument. When preforming
i.i.d.\ bond percolation with retention probability $p$
-- indicated by writing $\P_p$ for the probability measure --
on a rectangle $R$, the event
${\tt cross}_{\tt H}(R)$ is defined as the set of those subgraphs of the induced subgraph $\Z^2(R)$ that contains a
path between the two different vertical sides called a {\bf
horizontal crossing} in $R$. The event ${\tt cross}_{\tt V}(R)$, which we define by
interchanging the words vertical and horizontal above is called a {\bf
vertical crossing } in $R$.

If $R$ is any rectangle from ${\tt Rect}$, then let ${\tt ray}(R)$ be  the unique infinite path in $\T_{\tt RF}$ starting from $R$ (as a vertex of $\T_{\tt RF}$), the vertices in ${\tt ray}(R)$ form a sequence of rectangles (in an alternating way from ${\tt Rect}_{\tt V}$ and ${\tt Rect}_{\tt H}$) in which every consecutive pair is either well joined in a vertical to horizontal or in a horizontal to vertical way. Let ${\tt road}(R):=\bigcup_{Q \in {\tt ray}(R)}\Z^2[Q]$ be the graph union of the subgraphs induced by the rectangles in ${\tt ray}(R)$. 

We show that ${\tt RF}$ is robust by showing that already ${\tt road}(R)$ is robust for any $R\in {\tt Rect}$.

If we know a lower
bound for $\P_p({\tt cross_{\tt V}}(Q))$ for $Q\in V({\tt ray}(R))\cap {\tt Rect}_{\tt V}$ and $\P_p({\tt cross_{\tt H}}(Q))$ for $Q\in V({\tt ray}(R))\cap {\tt Rect}_{\tt H}$, then by the well-known Harris–FKG
inequality which states that for i.i.d. percolation any two increasing events are positively correlated
 (see, e.g.,
\cite[Thm.\ 2.4]{G}) we can get a lower bound for the probability of an infinite path simply
by multiplication, since if the appropriate (horizontal ones for horizontal rectangles, vertical ones for vertical rectangles) crossing events hold simultaneously, then from the individual paths within the rectangles we can extract the existence of a single infinite path (the situation is illustrated on the right of Figure~\ref{figChosenRoad}).
Such a lower bound is provided by the following lemma (stated as {Corollary 3.4} in \cite{OP}, and derived from an estimate by Bollobás-Riordan in \cite{BR}):

\begin{lemma}  \label{l.BRcorr}
For any fixed $p\in (\frac{1}{2},1)$, we can find
a $c>0$, a positive integer $\nu_0$ and a  $ \gamma
>0 $ such that for any positive integer $L$ the following holds.
If $R$ is an $Ln \times n$ rectangle where $ n > \nu_0$,
then $\P_{p}({\tt cross}_{\tt H}(R)) \geq c^{L/{n^{\gamma}}} $.
\end{lemma}

Of course a similar statement holds for $n \times Ln$ rectangles regarding the vertical crossings.

In \cite{OP} (section 3.6) the robustness of ${\tt road}(R)$ is shown by the above strategy (we assume that ${\tt j}(R)=0$, so that the indices below matches with $m_i$). That is, the parameters $(l_0,d_0),L_1,L_2,\dots$ can be chosen in such a way that if the $j$th rectangle in ${\tt road}(R)$ is a $\Lambda_j n_j \times n_j$ one (if it is horizontal) or a $n_j \times \Lambda_j n_j$ one (if it is vertical), then $\sum_{j=0}^{\infty}{\Lambda_j/{n_j^{\gamma}}}<\infty$, which makes the product $\prod_{j=0}^{\infty}c^{\Lambda_j/{n_j^{\gamma}}}$ to be positive. Since $\prod_{j=0}^{\infty}c^{\Lambda_j/{n_j^{\gamma}}}$ is a lower bound for the existence of an infinite path within ${\tt road}(R)$, by Harris-FKG in a Bernoulli $p$ percolation, ${\tt road}(R)$ is robust.

If we accept that those  calculations are correct in the case of ${\tt RF}$, then, when we move from ${\tt RF}$ to ${\Phi}_{\beta}({\tt RF})$, we will not need to know what $\Lambda_i$ is because we can use a soft argument: if $\sum_{i=0}^{\infty}a_i<\infty$, then there exists a sequence of positive integers $k_1\leq k_2 \leq \dots$, with $k_i\rightarrow\infty$ such that $\sum_{i=0}^{\infty}a_i k_i <\infty$ still holds. One construction is to select indices $j_1<j_2<\dots $ in such a way that $\sum_{j_i}^{\infty}\leq {1\over {2^i}}$, then let $k_i:=s$ if $j_s\leq i<j_{s+1}$. Then $\sum_{i=0}^{\infty}a_i k_i <\sum_{i=0}^{j_0}a_i+\sum_{i=0}^{\infty}{i\over {2^i}}<\infty$.

The reason this will complete our job is that ${\Phi}_{\beta}({\tt RF})$ still consists of folded slices which are isomorphic to planar rectangles, and the way these folded slices overlap is also isomorphic to the way how two well-joined planar rectangle overlap. Thus a crossing argument will also be applicable to the folded slices in ${\Phi}_{\beta}({\tt RF})$ as well. When we construct ${\Phi}_{\beta}({\tt RF})$ from ${\tt RF}$, then we have the sequence  $m_0\leq m_1\leq \dots$, with $m_i\rightarrow \infty$ as a parameter which we have great freedom to choose. 
When we cut a $\Lambda n \times n$ rectangle into $m$ horizontal slices, then a slice will approximately be a $\Lambda n \times n/m$ rectangle  (if $m$ does not divide $n$ it is only an approximation), and assuming that $n/m>\nu_0$ still holds, the lower bound from Lemma~\ref{l.BRcorr} for the horizontal crossing probability is $$c^{\Lambda /{({n\over {m}})^{\gamma}}}=c^{{m^{\gamma}\Lambda}\over{n^{\gamma}}}.$$

When we move from ${\tt RF}$ to $\Phi_{\beta}({\tt RF})$, then we do not only cut a $\Lambda n \times n$ rectangle into $m$ slices (which makes slices of $\Lambda n \times n/m$ rectangles) but we also fold those slices, and a folded slice is then isomorphic to a $\Lambda n+h \times n/m$ rectangle, where $h$ is some integer at most $3$ (corresponding to the need to moving up and down between the top and bottom layer of $\Z^2\times \{0,1\}$, when we fold the $\Lambda n \times n/m$ slice). However, if we ignore this additive term for a moment, we basically get that it is enough to choose the sequence  $m_0\leq m_1\leq \dots$, with $m_i\rightarrow \infty$ in such a way that  $$\sum_{j=0}^{\infty}{{m_j^{\gamma}{\Lambda_j}}\over {n_j^{\gamma}}}<\infty,$$ since this will imply by a crossing argument that the clusters of $\Phi_{\beta}({\tt RF})$ are robust.

Since $\sum_{j=0}^{\infty}{\Lambda_j/{n_j^{\gamma}}}<\infty$, we can find a sequence of positive integers $k_1\leq k_2 \leq \dots$, with $k_i\rightarrow \infty$, such that $\sum_{j=0}^{\infty}k_j{\Lambda_j/{n_j^{\gamma}}}<\infty$, so if we choose $m_j\rightarrow \infty$ in such a way that $m_j\leq k_j^{1/{\gamma}}$, then we are done, if we make sure that $n_j/m_j>\nu_0$, that is if $m_j<n_j/\nu_0$ (which is needed to apply Lemma~\ref{l.BRcorr}).
To ensure that $m_j<n_j/\nu_0$ holds, note that if the sequence $n_j$ is nondecreasing and $n_j\rightarrow \infty$, 
then with the above $k_j^{1/\gamma}$ sequence it will be true that ${\tt min}(n_j,k_j^{1/\gamma})$ is nondrecreasing and ${\tt min}(n_j,k_j^{1/\gamma})\rightarrow \infty$ still holds, so we can choose $m_j\leq {\tt min}(n_j,k_j^{1/\gamma})$ in  way that $m_0\leq m_1\leq \dots$ and $m_i\rightarrow \infty$ still holds. Finally, let us note that is indeed possible to choose the parameters $(l_0,d_0),L_1,L_2,\dots,$ in such a way that $n_j$ (which is the shorter side of the rectangles used) is strictly increasing.

Note also that so far we just used the approximate length of the sides of the slices as if they were $\Lambda n \times n/m$, when in reality both sides can be slightly different (partially by the extra length obtained from the folding and also from the possibility that $m$ does not divide $n$). Note that the crossing probabilities could only decrease if the longer side of the rectangle is lengthened and the shorter side is shortened, and then we can be sure that the longer side of the folded rectangle still will be less than $2\Lambda n$ and the shorter side cannot be shorter than $n/2m$, thus if $n/2m>\nu_0$ still holds, then we can use Lemma~\ref{l.BRcorr} to estimate the crossing probability of a $2\Lambda n \times n/2m$ rectangle (and use it as a lower bound for the crossing probability of the folded slice). In this way the sum whose convergence will ensure the robustness of our clusters will only increase by a constant multiplier compared to the estimates we used before (and the condition $n/2m>\nu_0$ is also easy to satisfy). This completes the proof of Theorem~\ref{t.slab}.

\section{The case of $\Z^2$}\label{s.planarNO}

In this section we prove Theorem~\ref{t.plane}.
We will use some basic tools from the theory of invariant processes in groups in general (for more background see the research paper \cite{BLPS} or the the books \cite{LPbook}, \cite{PGG}).
The {\it number of ends} in a graph $G$ is the supremum of the number of infinite components of $G\setminus K$ over all finite subsets $K$ of $G$ (for an explanation what is an end in general see section 7.3 in \cite{LPbook}).
If the number of ends in a graph is $k$, we will call it ``$k$-ended" (for trees it is consistent with the notion of one-ended trees we used before).
 The following is a well-known fact (stated as an exercise in \cite{PGG}, following the landmark Burton-Keane argument):

In an invariant percolation on an amenable group (such as $\Z^2$) there cannot be a cluster with more than two ends.

We will also use the Mass-Transport Principle for the case of countable groups (see for example \cite{LPbook}). We only need the following consequence: an invariant process on a countable group cannot produce any infinite cluster in which a non-empty finite set of vertices are distinguished. This is true, because in every infinite cluster $C$ where a non-empty finite set of vertices is distinguished we could pick a single vertex $v_C$  uniformly and every vertex from $C$ could send mass $1$ to $v_C$ in which case the expected mass sent out is at most $1$ but the expected mass received is infinite.

Let the {\it dual lattice} of $\Z^2$ be the graph (isomorphic to $\Z^2$) whose vertices are of the form $a+(1/2,1/2)$ for $a\in \Z^2$, and edges are $\{a+(1/2,1/2),b+(1/,1/2\}$ for $\{a,b\}\in \Z^2$. There is a natural bijection between the edges of the two graphs, namely the ones which are crossing each other can be put into pairs (that is $\{a,a+1\}$ for $a\in \Z$ is the pair of $\{a+1/2,a-1/2\}$ while $\{a,a+1\}$ is the pair of $\{a-1/2,a+1/2\}$).
We will call configurations $\omega \in \{0,1\}^{E(\Z^2)}$ in the ``original" $\Z^2$ {\it primal} and we associate to such an $\omega$ a natural {\it dual} $\omega^*$: let an edge $f$ in the dual graph be in $\omega^*$ iff its pair $e$ is not in $\omega$. 

If $C_1$ and $C_2$ are two disjoint infinite components of $\Z^2$, then we can find an at least two ended dual component between them (thus when it arises from an invariant process, it is exactly two-ended). This is because if $P$ is a shortest path from $C_1$ to $C_2$, and we start an exploration of the boundary of $C_1$ from the first edge of $P$ into both direction, then this exploration finds an infinite path in the dual on either side of $P$, and if we delete the dual edges corresponding to the edges of $P$, then these two infinite paths end up in different components in the dual, thus the dual region between $C_1$ and $C_2$ must be at least two-ended.

Let us say that a dual cluster $C_1$ {\it touches} a primal cluster $C_2$ if $C_1$ has a vetex $v_1$ for which $v_1=v_2+(\sigma_1/2,\sigma_2/2)$ for a vertex $v_2$ of $C_2$, where $\sigma_i\in \{-1,+1\}$.

Assume that ${\cal I}$ is an invariant percolation process on $\Z^2$ which produces at least one robust cluster. We will show that it cannot have more than two robust clusters. The argument is cleaner if we assume that ${\cal I}$ has robust clusters and isolated vertices only. So let us delete the edges from any non-robust cluster. 
In what follows ${\cal I}^*$ is the dual process (which produces a dual configuration of ${\cal I}$).

\begin{figure}[htbp]
\centerline{
\includegraphics[width=0.5\textwidth]{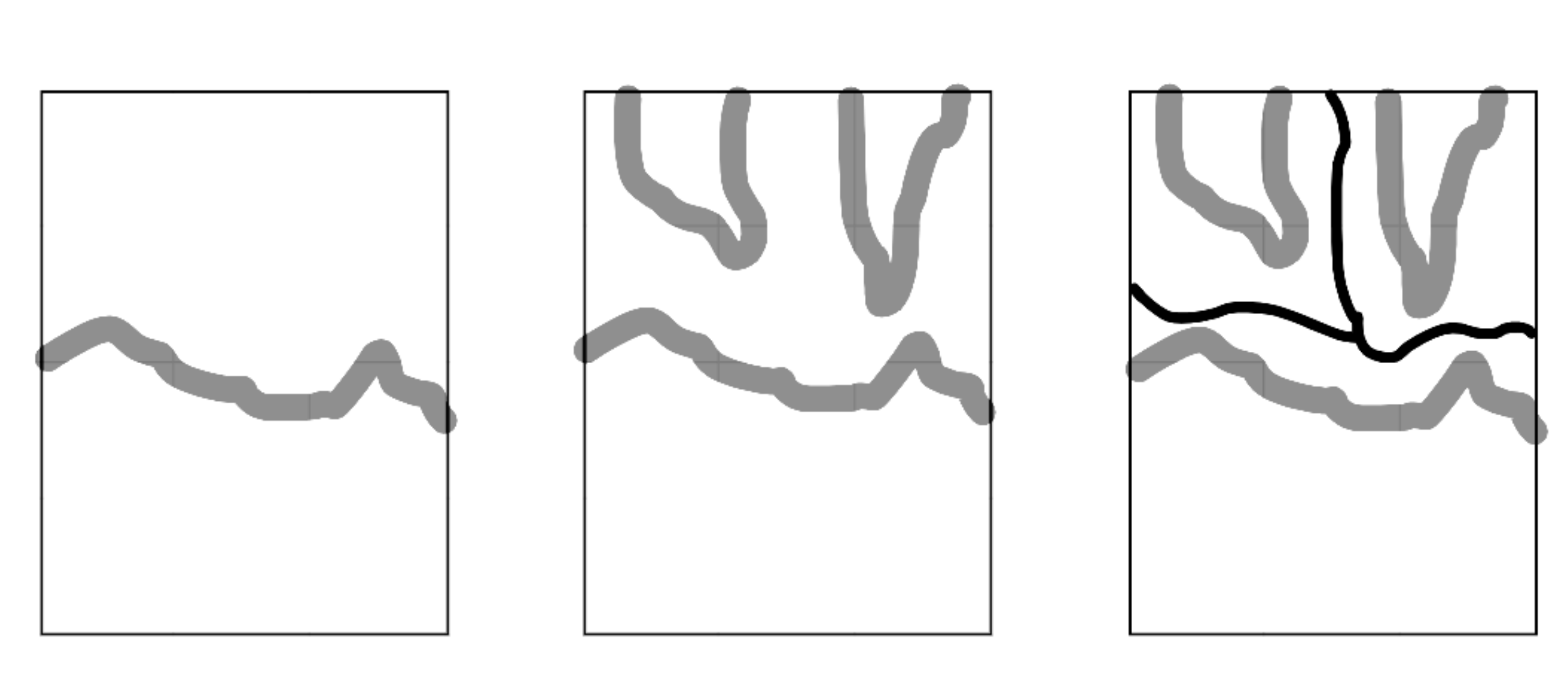}
}
\caption{The picture shows the following steps: if there are at least two infinite primal cluster, then there must exists at least one 2-ended dual clusters ${\cal D}$ which separates the primal $\Z^2$ into two infinite pieces. If they were more than one infinite dual clusters neighboring ${\cal D}$ on the same side of ${\cal D}$, then there would exists a primal cluster with more than two ends contradicting the amenability of $\Z^2$.}
\label{planarTHREEend}
\end{figure}

If ${\cal I}$ has at least two infinite clusters, then there must exist a 2-ended cluster of ${\cal I}^*$ which separates them. Let ${\cal D}$ be any 2-ended cluster of ${\cal I}^*$, this separates the primal $\Z^2$ into two infinite pieces (and possibly into finite ones). Let us call two dual clusters {\it neighbors} if there is an infinite primal cluster that they both touch and consider the infinite dual clusters neighboring ${\cal D}$. By inspecting the number of ends in the primal configuration (see Figure~\ref{planarTHREEend}) there can only be at most one such neighbor on either side of ${\cal D}$ (otherwise the primal configuration must contain components with at least three ends). 

Call them ${\cal N}_1({\cal D})$ and ${\cal N}_2({\cal D})$, if they exist. We will see that, using the robustness of the clusters of ${\cal I}$, not even a single such neighbor can exist. 
Assume that at least one of them exists, let ${\cal N}_{\cal D}$ be one of them selected uniformly at random, then the primal component ${\cal C}_{\cal D}$ between ${\cal D}$ and ${\cal N}_{\cal D}$ must be two ended. 

Let the minimal distance between ${\cal D}$ and ${\cal N}$ be $m_{\cal D}$. Let $M({\cal D})$ be the set of those vertices of ${\cal D}$ for which there exists a vertex $y$ of $N$ such that ${\tt dist}(x,y)=m_{\cal D}$. Because the association of $M({\cal D})$ to two-ended dual components ${\cal D}$ is an invariant selection of a non-empty (if at least one of ${\cal N}_i({\cal D})$ exists) set of vertices, $M({\cal D})$ must be infinite.

This implies that the primal cluster ${\cal C}_{\cal D}$ between ${\cal D}$ and ${\cal N}_{\cal D}$ cannot be robust, contradicting the assumption on ${\cal I}$. This is because if $P_{\tt min}({\cal D},{\cal N}_{\cal D})$ is the collection of the infinitely many shortest dual paths between ${\cal D}$ and ${\cal N}_{\cal D}$ and ${\tt inv}(P_{\tt min}({\cal D},{\cal N}_{\cal D}))$ is any invariantly selected non-empty subset of $P_{\tt min}({\cal D},{\cal N}_{\cal D})$, then deleting all the edges of ${\cal C}$ corresponding to the edges of the paths in ${\tt inv}(P_{\tt min}({\cal D},{\cal N}_{\cal D}))$ cuts ${\cal C}$ into finite components and a Bernoulli $p$ percolation with $p<1$ on ${\cal C}$ will indeed delete such collection of edges from ${\cal C}$.

\section{Acknowledgements} This work was partially supported by the ERC Consolidator Grant 772466 ``NOISE''.
I am indebted to Gábor Pete for useful discussions throughout, and to Russell Lyons to ask the question answered by Theorem~\ref{t.plane}.

\end{document}